\documentclass[10pt,a4paper,reqno]{amsart}
 \usepackage{amscd}
 \usepackage{amsmath}
  \usepackage{amsthm}
 \usepackage{epsf}
 \usepackage{epsfig}
 \usepackage{graphics}
 \usepackage{graphicx}
 \usepackage{geometry}
 \usepackage{tikz-cd}
 \usepackage{tikz} 
 \usepackage{amssymb}
 \usepackage{textcomp}
 \usepackage{gensymb}
 \usepackage{tkz-euclide}
 \usepackage{verbatim}
 \usepackage[english]{babel}
 \usepackage{tabu}
 \usepackage{amsmath}
 \usepackage{hyperref}
 \usepackage{cleveref}
 \usepackage{setspace}
 \usepackage{caption}
\usepackage[verbose]{placeins}

 \usepackage[T1]{fontenc}
\usepackage{array}
\usepackage{booktabs}

\usepackage{hyperref}
\hypersetup{
	colorlinks,
	citecolor=red,
	filecolor=black,
	linkcolor=blue,
	urlcolor=black
}

 \numberwithin{equation}{section}
 \newcommand{\field}[1]{\ensuremath{\mathbb{#1}}}
 
 \newcommand{\PP}{\field{P}}
 
 \newcommand{\RR}{\field{R}}
 \newcommand{\ZZ}{\field{Z}}
   \newcommand{\M}{\mathcal{M}}
   \newcommand{\Z}{\mathcal{Z}}
    
   \newcommand{\PL}{\mathcal{P}}

 \newcommand{\lra}{\longrightarrow}
 
 \providecommand \hf{\hspace*{0.5cm}}
 
 \theoremstyle{plain}
 \newtheorem{thm}{Theorem}[section]
 
 \newtheorem{lmm}[thm]{Lemma}
 \newtheorem{prop}[thm]{Proposition}
 \newtheorem{que}[thm]{Question}
 \newtheorem{eg}[thm]{Example}
 \newtheorem{rem}[thm]{Remark}

 \newtheorem{defn}[thm]{Definition}

 \usepackage{hyperref}
 \usepackage{float}
 \usepackage{array,tabularx}
 \usepackage{babel,tabularx,ragged2e,booktabs}

\begin{document}
 \title{Tropical rational curves 
with first order tangency via WDVV}

 \author[Anantadulal Paul]{Anantadulal Paul}
\address{School of Mathematics,
Tata Institute of Fundamental Research,
Dr. Homi Bhabha Road,
Navy Nagar, Mumbai 400 005,
Maharashtra, India.}
\email{paul@math.tifr.res.in/ paulanantadulal@gmail.com}
\author[Aditya Subramaniam]{Aditya Subramaniam}
\address{Department of Mathematics, Manipal Institute of Technology, Manipal Academy of Higher Education, Manipal, Karnataka 576104, India.}
\email{aditya.subramaniam@manipal.edu / adisubbu92@gmail.com}




 \begin{abstract} 
In this article, we study the tropical counterpart of the enumeration of rational curves in $\mathbb{CP}^2$ with first order tangency. We use the tropical analogue of the WDVV technique to compute rational tropical plane curves of degree $d$ tangent to a degree $l$ tropical plane curve and passing through $3d-2$ points in general position. As Mikhalkin's correspondence  (cf. \cite{Mik-g=0-char-numbers})  suggests, our numbers agree with earlier results on tangency in complex geometry. 
 \end{abstract}
\maketitle

\section{Introduction}
The number of complex plane curves of any genus  tangent to an ample divisor and passing through an appropriate number of generic points is known as characteristic numbers. In complex projective plane $ \mathbb{CP}^2$, the enumeration of characteristic numbers is a classical problem. In modern techniques, this theory advances by the name relative Gromov-Witten theory. After Mikhalkin's seminal paper \cite{Mik-correspondence} and his famous so-called ``Correspondence Theorem," tropical geometry gets a lot of attraction from various fields of mathematics. For instance, a remarkable result in tropical geometry is the tropical reformulation of the classical enumeration of  degree $d$  plane curves of genus $g$ passing through $3d-1+g$ points in general position (see \cite{Mik-correspondence}).\\
\hf In recent times, there has been considerable research in understanding the tropical counterpart of several enumerative questions involving tangencies (contact orders). In the present paper, we consider the tropical counterpart of the following classical question in complex geometry :
\begin{que}
\label{First-ord-tang-que}
   What is the number of rational plane curves of degree $d$ passing
through $3d-2$ generic points and tangent to
a degree $l$ curve?
\end{que}

On the complex geometric side, Caporaso and Harris  defined the concept of relative Gromov-Witten invariants and completely solved the general question of 
counting plane curves of a given genus and degree with
local contact orders to a fixed line 
satisfying appropriate point constraints
(see \cite{CH}). Later, in his papers \cite{AnGaPhD}, \cite{AnGa1} and \cite{AnGa2},  Gathmann gives a systematic approach for enumeration of characteristic numbers of rational and elliptic degree $d$ curves in 
$\mathbb{CP}^n$ with multiple tangencies to an ample hypersurface $Y \subset \mathbb{CP}^n$ at multiple points of various orders. 
As a remarkable application, he used this approach to enumerate Gromov-Witten invariants of the quintic threefold.\\ 
\hf On account of these classical results, it is natural to ask what would be the tropical counterpart of these. Several instances of this in the tropical side have already appeared. For instance, Gathmann and Markwig (see \cite{Trop-CH}) studied the tropical analogue of the classical Caporaso-Harris result. Very recently, there have been many advances appeared in the literature on these type of problems. Amongst others, Markwig et al. studied the enumeration of tropical bitangents, and along the way,  showed all possible deformation classes of tropical bitangents for quartics. They also considered the $\mathbb{A}^1$-homotopy theoretic invariants corresponding to the tangencies (see the results \cite{Lift-trop-bitangent}, \cite{Trop-quartic-bitangent}). \\
\hf In \cite{Mik-g=0-char-numbers}, Mikhalkin et al. established a special case of the correspondence theorem that ensures
that the Question \ref{First-ord-tang-que}
is a well posed problem in tropical geometry i.e., the classical characteristic numbers in complex geometry coincide with the tropical enumeration. They used the floor
diagram calculus to reduce the problem into pure combinatorics.
Question \ref{First-ord-tang-que} has already been studied using complex geometric techniques due to, amongst others, 
\cite{AnGaPhD}, \cite{Pan@Qdiv}, and in the tropical side by \cite{Mik-g=0-char-numbers}.\\
\hf \hf In this paper, we solve Question \ref{First-ord-tang-que} tropically by using the tropical analogue of the WDVV equation, closely following \cite{GathTropical}. Our method is extremely straightforward and geometric. We hope that our method can be generalized to retrieve the complete result in \cite{Mik-g=0-char-numbers}. This will be explored elsewhere. 

\section{Overview of our result}
In this paper, the main technical idea that we use to study Question \ref{First-ord-tang-que} tropically is to adapt  the classical approach of enumerating complex plane rational curves using the WDVV technique.

\hf We now recall the classical degeneration process which automatically appears while considering the WDVV technique to enumerate
 rational curves in 
$\mathbb{CP}^2$ with first order tangency to a fixed degree $l$ curve.
The situation can be described by the following picture:
\begin{figure}[hbt!]
\begin{center}\includegraphics[scale = 0.6]{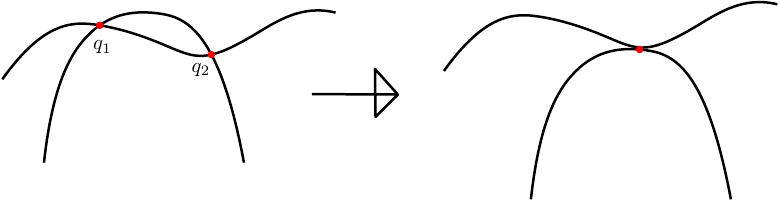}
\end{center}
\caption{Two distinct transverse intersections limiting to a first order tangency }
\label{pic idea}
\end{figure}
\FloatBarrier
Denote a degree $l$ curve  inside $\mathbb{CP}^2$ by $\mathcal{L}^{l}$.
The right hand side (RHS, in short) in the above picture shows a rational plane curve tangent to $\mathcal{L}^{l}$, while the left hand side (LHS, in short) represents  a rational curve intersecting $\mathcal{L}^{l}$ at two distinct points transversally. The picture depicts  a typical rational degree $d$ curve intersecting  $\mathcal{L}^{l}$ transversally at two distinct points $q_1, q_2$, limiting to a typical rational degree $d$ curve tangent to $\mathcal{L}^{l}$
at some point. Note that the intersection multiplicity of the two curves at this point is two. In the case when $l=1$, the above idea has been explored in \cite{IRAA-arX}.  

\subsection*{WDVV algorithm}
 Let us first describe the main idea of deriving the recursive formula to compute the numbers for rational degree $d$ curves tangent to a fixed degree $l$ curve. Let us recall that $\overline{M}_{0,0}(\mathbb{CP}^2, d) $ be the Deligne-Mumford compactification of the Kontsevich moduli space of rational curves (with no marked points). Let $\mathcal{H}$ denote the divisor that corresponds to the subspace of curves that pass through a generic point. We note that the intersection number defined by
\begin{align*}
n_d:= deg\left(\big[ \overline{M}_{0,0}(\mathbb{CP}^2, d)  \big] \cdot \mathcal{H}^{3d-1}\right),
\end{align*}
computes rational degree $d$ curves through $3d-1$ general points. The numbers $n_d$ are given by an explicit recursion formula due to the famous work due to Kontsevich and Manin \cite{K.M} and independently by Ruan and Tian \cite{RT}. The numbers $n_d$ play the role of the base case for our desired recursion formula.


\hf Let $\mathsf{X}$ be the subspace of ${M}_{0,n}(\mathbb{CP}^2, d) $ consisting of rational degree $d$ curves (with smooth domains) where the image of the first two marked points lie on two distinct points of a fixed degree $l$ curve $\mathcal{L}^l$ i.e., the rational curve intersects $\mathcal{L}^l$ transversely at two distinct points. Denote by $\overline{\mathsf{X}}$ the closure of $\mathsf{X}$ inside $ \overline{M}_{0,n}(\mathbb{CP}^2, d)$. 
Next, we will be applying the WDVV technique.

\hf \hf Consider $\overline{M}_{0,4}$, the moduli space of four pointed stable rational curves. This is isomorphic to $\mathbb{P}^1$, hence path connected. Therefore, any two points in this space determine the same divisor. Let $D(ij \mid kl)$ be the divisor in $\overline{M}_{0,4}$ representing the two componented rational curves where the marked points $(x_i,x_j)$ lie in one component and $(x_k, x_l)$ lie on the other component. We occasionally denote this by $(ij \mid kl)$. Thus, we have the following divisorial identity in $H^{\ast}( \overline{M}_{0,4}, ~\ZZ)$. Thus
\[ \mathcal{R}:=  D(i j | k l)  -  D(i k | j l) \equiv 0 \] as divisors. 
Let us consider the forgetful map 
\begin{align*}
\pi: \overline{M}_{0,n}(\mathbb{CP}^2, d) &\longrightarrow \overline{M}_{0,4}.
\end{align*}
We shall remind the reader that this is not exactly a forgetful map; this forgets the map, and if the resulting domain
is unstable, the map stabilizes that.
Then
\begin{align}
  \pi^{\ast} D(i j | k l) \equiv \pi^{\ast} D(i k | j l).  
\end{align}
Let us define a class $\Z$ in $ \overline{M}_{0,n}(\mathbb{CP}^2, d) $ as 
\begin{align*}
    \Z:=  ev_3^{\ast}([pt]) \cdot ev_4^{\ast}([pt]) \cdot \big\{ ev_5^{\ast}([pt]) \cdots ev_n^{\ast}([pt])\big\}.
\end{align*}
We now consider the intersection theory of the cycle $\Z$, intersected with 
\[  \pi^{\ast} \mathcal{R} \cdot \overline{ \mathsf{X} } \] 
inside the moduli space $\overline{M}_{0,n}(\mathbb{CP}^2, d)$.

Hence, 
we get
\begin{align}
\big[\overline{ \mathsf{X} } \big]\cdot \big[\pi^*(12|34)\big] \cdot \Z  & = 
\big[\overline{\mathsf{X}}\big]\cdot \big[\pi^*(13|24)\big] \cdot \Z. \label{T1_WDVV}
\end{align}
Observe that this is an equality of numbers. Next, by unwinding the LHS and RHS of \cref{T1_WDVV}  
we will get our desired recursion formula.\\
\hf In this paper, we employ the above idea in tropical geometry modelled on the influential work by Gathmann and Markwig \cite{GathTropical}. Analogous to the classical moduli space $\overline{M}_{0,n}(\mathbb{CP}^2, d)$, they introduced the moduli space of stable tropical degree $d$ plane curves with $n$ markings $\M_{d,n}$. The authors showed that $\M_{d,n}$ with $n \geq 4$ admit forgetful maps to $\M_4$ (the tropical analogue of $\overline{M}_{0,4}$). Then, the WDVV approach in the tropical setting reduces to a degree computation of a suitably defined tropical map (a morphism between polyhedral complexes of pure dimension) incorporating all the required incidence conditions at two different generically chosen points with two different types of $\M_4$ coordinates having large intrinsic lengths.\\
\hf Our main result (see \Cref{final theo}) follows from the treatment taken in \cite{GathTropical} by modifying their setup appropriately. More precisely, we define a new tropical map $\tilde{\pi}$ (see Equation \eqref{map-giving-numbers}) whose degree at two generically chosen points in $\mathbb{R}^{2n-4}\times \M_{4}$ produces a recursion formula for first order tangency. We will follow their tropical presentation very closely. The idea behind tropical rederivation of the famous Kontsevich's recursion formula in \cite{GathTropical} is parallelly used to derive the Kontsevich type recursion formula for $N_d^{\mathsf{T}_1}(l)$ - The number of rational tropical degree $d$ plane curves tangent to a tropical curve of degree $l$ satisfying point constraints (see \Cref{KoTa}).

\section{Plane tropical curves}
In this section, we recall the notion of abstract tropical curves, plane tropical curves and their moduli spaces. We mainly follow the exposition and notations from \cite{GathTropical}. Recall if $\{ I_j: ~j= 1,\ldots,n \}$ be a set of closed, bounded, or half-bounded finite intervals, then a tropical graph $\Gamma$ is formed by identifying the end points of these intervals. An edge of $\Gamma$ will be called bounded or unbounded if its corresponding open interval is so. We will denote $\Gamma_{0}^{1}$ for a bounded edge and $\Gamma_{0}^{\infty}$ for an unbounded edge.
\begin{defn}
   An abstract tropical curve is a connected
graph $\Gamma$ of genus zero whose all vertices have valence at least $3$. An $n$ marked tropical curve is
a tuple $(\Gamma, x_1,\ldots,x_n)$ where $x_1,\ldots,x_n \in \Gamma_{0}^{\infty}$. 
\end{defn}
 Note that two such curves $(\Gamma, x_1,\ldots,x_n)$ and $(\Gamma^{\prime}, y_1,\ldots,y_n)$ are isomorphic if there is a homomorphism $\Gamma \longrightarrow \Gamma^{\prime}$ that sends $x_i$ to $y_i$ so that every edge of $\Gamma$ is mapped bijectively onto an edge of $\Gamma^{\prime}$ by
an affine map of the form $t \longrightarrow a \pm t$.  The moduli space of $n$-marked tropical curves will be denoted by $\mathcal{M}_n$ (Analogous to $\overline{M}_{0,n}$ in the classical setting). Notice that an element of $\mathcal{M}_4$ has the following four
possible combinatorial types:
\begin{figure}[hbt!]
\vspace*{0.1cm}
\begin{center}\includegraphics[scale = 0.7]{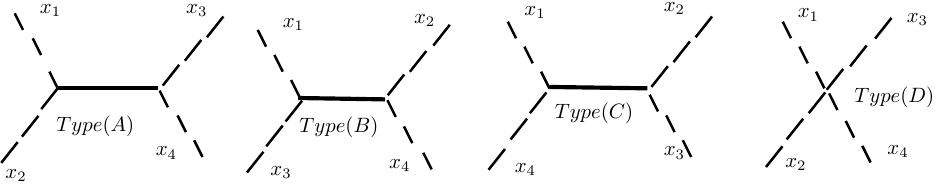}
\vspace*{-0.2cm}
\end{center}
\caption{Tropical curves of types (A), (B), (C), and (D)}
\label{type of curves in M4}
\end{figure}
\FloatBarrier

We now define the notion of an $n$-marked plane tropical curve.
\begin{defn}
   A $n$-marked plane tropical curve is a tuple $(\Gamma, x_1,\ldots,x_n, h)$,
where $x_1,\ldots,x_n \in \Gamma_{0}^{\infty}$ are distinct unbounded edges of $\Gamma$, and $h: \Gamma \longrightarrow \mathbb{R}^2$ is a continuous map, such that:
\begin{itemize}
    \item Restricted to each edge of $\Gamma$ the map $h$ takes the form $h(t) = a + t \cdot v$ for some $a \in \mathbb{R}^2$ and $v \in \mathbb{Z}^2$.
    \item  For every vertex $V$ of $\Gamma$, we have the balancing condition $$ \sum_{F: ~flags}{v(F)} = 0.$$
    \item  Unbounded edges $x_i$ gets mapped to a point in $\mathbb{R}^2$ by $h$.
\end{itemize}
\end{defn}
Note that two tropical plane curves are equivalent if the underlying abstract tropical curves are so. The degree of a plane tropical curve with $n$ markings is defined as the multiset $\bigtriangleup$ of directions of its non-marked unbounded edges. Denote the space of all plane tropical curves of degree $\bigtriangleup$ and $n$ markings by $\mathcal{M}_{ \bigtriangleup, n}$ (analogous to $\overline{M}_{0,n}(\PP^2, d)$ in classical setting). For more details and several examples, we refer the reader to \cite{GathTropical}.\\
\hf The combinatorial type of a marked abstract tropical curve
$(\Gamma, x_1,\ldots,x_n)$ is defined to be the homeomorphism class of $\Gamma$ relative to the markings $x_1,\ldots,x_n$. Similarly, the combinatorial type of a plane tropical curve $(\Gamma, x_1,\ldots,x_n, h)$ is the combinatorial type of
the marked abstract tropical curve together with the direction vectors $v(F)$ for all flags $F \in \Gamma^{\prime}$. Here $\Gamma^{\prime}$ denotes the set of all flags of $\Gamma$. The codimension of such a combinatorial type $\alpha$ in both cases is defined to be
\begin{align*}
    \text{Codim}~\alpha := \sum_{V\in \Gamma^{0}} \left(Val(V) - 3 \right),
\end{align*}
where $\Gamma^{0}$ denotes the set of vertices of $\Gamma$. Hence, we denote $\M^{\alpha}_n$ (respectively  $\M^{\alpha}_{\bigtriangleup, n}$) to be the subset of $\M_n$ (respectively subset of $\M_{\bigtriangleup, n}$) corresponding to marked tropical curves of type $\alpha$. 
\begin{prop}
    For any positive integer $n$ and combinatorial type $\alpha$, the
space $\M^{\alpha}_n$ (respectively $\M^{\alpha}_{\bigtriangleup, n}$) is naturally an (unbounded) open convex polyhedron in a real
vector space of expected dimension
\begin{align*}
    \text{dim}~\M^{\alpha}_n &= ~n-3- codim(\alpha)~~\\
     (\text{respectively dim}~\M^{\alpha}_{\bigtriangleup, n} &= ~\mid \bigtriangleup\mid -1 + n -codim(\alpha).)
\end{align*}
\end{prop}

\begin{proof}
    See \cite[Proposition 2.11]{GathTropical}.
\end{proof}
\hf As in the classical setting, we have the following forgetful morphisms 
\begin{align}
    \delta:  \M_{\bigtriangleup, n} \lra  \M_{m} , ~\text{and}~   \tilde{\delta}:  \M_{\bigtriangleup, n} \lra  \M_{\bigtriangleup, m}.
\end{align}
We refer the reader to \cite[Definition 4.1]{GathTropical} for more details on these two kinds of forgetful maps. \\ 
\hf We now define the notions of tropical multiplicity and degree of a morphism $f:X \lra Y$ between polyhedral complexes $X$ and $Y$ of equal pure dimension. For a general point $p\in X$, we define the multiplicity $mult_{p}(f)$ of $f$ at $p$ to be the absolute value of the determinant of the underlying linear map.
\begin{defn}
Let $f:X \lra Y$ be a morphism between polyhedral complexes of equal pure dimension. A point $Q\in Y$ is said to be in $f$-general position if $Q$ is in general position in $Y$ and the points in $f^{-1}(Q)$ are also generic in $X$.
\end{defn}
Next, we define the degree of $f$ at a $f$-general point $Q$ as
\begin{align}\label{deg-mult}
    deg_{Q}(f)  := \sum_{P \in f^{-1}(Q)}mult_{P}(f).
\end{align}
Note that in the case when both $X$ and $Y$ are the moduli spaces of abstract or plane tropical curves, the multiplicity and degree do not depend on the choices of coordinates on the cells of $X$ and $Y$ (see \cite[Remark 3.2]{GathTropical}).\\
\hf As customary, we have evaluation morphisms
\begin{align}
    ev_{i}: \M_{\bigtriangleup, n} \lra \RR^2, ~\hspace{1cm}~~ (\Gamma, x_1,\ldots,x_n, h) \mapsto h(x_i).
\end{align}
Note that these are morphisms of polyhedral complexes. Let us denote by $ev$ to be the product of $n$-evaluations $ev_1 \times \ldots \times ev_n: \M_{\bigtriangleup, n} \lra \RR^{2n}$. Let $n= \mid \bigtriangleup \mid -1$. Then the source of the map $ev$ and the target have the same dimension $2n$. Therefore, we can define the numbers
\begin{align}
    N_{\bigtriangleup}(w) := deg_{w}(ev)
\end{align}
for all points $w \in \RR^{2n}$ in $ev$-general position. These numbers are counting tropical curves of degree $\bigtriangleup$ passing through points $w$, where each curve $C$ is counted with
a certain multiplicity $mult_{C}(ev)$. Note that  \cite[Proposition 3.8]{GathTropical} shows how to compute the tropical multiplicity $mult_{C}(ev)$. \\

\begin{eg}
\label{trop_mult_eg}
Let us recall the tropical multiplicity calculation for the following combinatorial type $\alpha$:
\begin{figure}[hbt!]
\vspace*{0.1cm}
\begin{center}\includegraphics[scale = 0.7]{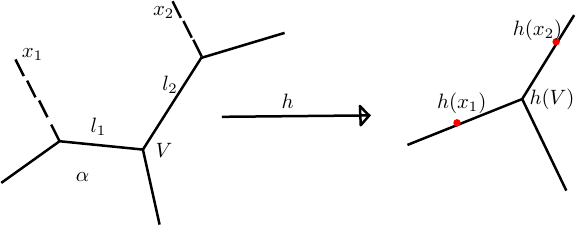}
\vspace*{-0.2cm}
\end{center}
\end{figure}
\FloatBarrier
Choosing the root vertex $V$, first fix the image coordinate $h(V) = (a, b)$. Then, from the root vertex, there are two bounded edges having lengths $l_1, l_2$ with their direction vectors $v_1 = (v_{1,1}, v_{1,2})$ and $v_2 = (v_{2,1}, v_{2,2})$. Thus, $a, b, l_1, l_2$ are the coordinates for the type $\alpha$ of $\M_{\Delta, 2}^{\alpha}$. This gives us the following matrix for the map $ev = ev_1 \times ev_2$ where $h(x_i) = h(V) + l_i \cdot v_{i},~~ i = 1, 2$\\
\[\begin{pmatrix}
1 & 0 & v_{1,1} & 0\\
0 & 1 & v_{1,2} & 0\\
1 & 0 & 0 & v_{2,1}\\
0 & 1 & 0 & v_{2,2}
\end{pmatrix}\]
Hence, the multiplicity of the $ev$ map is given by the absolute value of the determinant of this matrix, i.e. $mult_{\alpha}(ev) = |det(v_1, v_2)|$.
\end{eg}


\section{Intermediate intersection numbers}
We now consider the variant of constraints where the rational tropical plane curve passes through general points with one of the mark points evaluating to a degree $l$ curve and two of the mark points evaluating to two distinct points in a degree $l$ curve. We will make an abuse of notation and denote the fixed degree $l$ plane tropical curve by $\mathcal{L}^l$. Let $N_d^{\mathsf{T}_0}$ be the number of rational degree $d$ tropical plane curves intersecting $\mathcal{L}^l$ transversally at some point and passing through $3d-1$ generic points. This enumerative invariant amounts to counting tropical curves passing through $3d-1$ generic points with some of the mark points lying on a one-dimensional cycle inside $\mathbb{R}^2$. For any nonnegative integer $n$, we define $N_d^{\mathsf{T}_0}(n)$ to be the number of rational degree $d$ tropical plane curves that intersect $\mathcal{L}^l$ transversally at some point, where the fixed degree $l$ curve intersects a class  $\mathsf{L}^n$ and passes through $3d-1-n$ points in general position. Here, $\mathsf{L}$ denotes the class of a line inside tropical plane. Hence, $N_d^{\mathsf{T}_0}(0) = N_d^{\mathsf{T}_0}$. When $n = 1$, we often denote the number $N_d^{\mathsf{T}_0}(1)$ by $N_d^{\left(  \mathsf{T}_0 \right)_{\textnormal{pt}}}$. This number counts rational degree $d$ tropical plane curves intersecting $\mathcal{L}^l$ at a fixed point (where the fixed point is seen as a point of intersection of 
two general line classes) and passing through additional $3d-2$ points in general position.

\hf By abuse of notation, $n_d$ denotes the number of plane tropical rational degree $d$ curves through $3d-1$ general points. Recall that these numbers are enumerated in \cite{GathTropical} via the degree computation by considering the map 
\begin{align}
\label{TKon_num}
ev = ev_1 \times \ldots \times ev_n : \M_{\bigtriangleup, n} \lra  \RR^{2n}.
\end{align}
Hence, $n_d = deg_{C}(ev)$, where $C$ is a curve in $ev$-general position (cf. \cite[Section 5]{GMK09}).\\

Note that the intersection number $N_d^{\mathsf{T}_0}(n)$ can be seen as the degree of the following map
\begin{align}
\phi_{n}:= \left(ev \times ev_{n+1} \times \delta \right) = ev_1\times ev_2\times ev_3 \times \ldots \times ev_{n+1} \times \delta: \M_{\bigtriangleup, n+1} \lra  \RR^{2n} \times \RR \times \M_{4},
\end{align}
where the first $n$ evaluations are mapped to $n$ general points in $\mathbb{R}^2 $ with the $(n+1)^{th}$ marking mapping to $\mathsf{L}^n$. Thus $N_d^{\mathsf{T}_0}(n) = deg(\phi_{n} )$.
\begin{lmm}\label{ndt0}
For $n \geq 0$, 
$N_d^{\mathsf{T}_0}(n)$ is given by
\begin{align}
\label{T0_num}
N_d^{\mathsf{T}_0}(n) = \begin{cases}
l d  ~n_d, ~~\qquad \qquad &\text{for}~ n = 0\\
n_d, ~~\qquad \qquad &\text{for}~ n = 1\\
0,   ~~~~~\qquad \qquad \qquad &\text{for}~ n \geq 2.
\end{cases}
\end{align} 
\end{lmm}

\begin{proof}
Recall that the number $n_d$ is defined as $deg_{C}(ev)$, where $C = (\Gamma, x_1, \ldots, x_{n}, h)$ is a generic plane tropical curve and the map $ev$ is defined as in \Cref{TKon_num} (see \cite[Theorem 5.6]{GathTropical} for more details). Next, consider a generic plane tropical curve $C = (\Gamma, x_1, \ldots, x_{n+1}, h)$ with the extra $l ~\mathsf{L}$ class as incidence condition, i.e.,  when the marking $x_{n +1}$ is mapped to the curve $\mathcal{L}^l$ intersecting the class $\mathsf{L}^n$.
Note that $x_{n+1}$ has to be adjacent to a  trivalent vertex (cf. \cite[Lemma 3.6]{Hanah_decendent}). Therefore, applying the result \cite[Lemma 6.6, and Theorem 5.3]{Hanah_decendent}, we have $$deg_C(\phi_{n}) = \left( \mathsf{L}^n \cdot \tilde{C} \right)_{h(x_{n+1})} \cdot deg_{\tilde{C}}(ev ),$$ where the curve $\tilde{C}$ is obtained from $C$ by forgetting $x_{n+1}$ ( i.e., by straightening the two valent vertex while removing $x_{n+1}$ from $C$), and the map $ev$ records the remaining $n$ point evaluations. Thus, any curve $C$ appearing in the consideration of $deg_C(\phi_{n})$ leads to a curve $\tilde{C}$ by forgetting the marked end $x_{n+1}$. Conversely, one can choose a point $x \in \left( \mathsf{L}^n \cdot \tilde{C} \right) $ and by attaching a marked end, it is possible to get $C$ appearing in $deg_C(\phi_{n})$. Using tropical B\'{e}zout's Theorem (see \cite[Theorem 1.3.2]{Strum-book}), we see that 

\begin{align*}
\left( \mathsf{L}^n \cdot \tilde{C} \right) = \begin{cases}
l d , ~~\qquad \qquad &\text{for}~ n = 0\\
1, ~~\qquad \qquad &\text{for}~ n = 1\\
0,   ~~~~~\qquad \qquad \qquad &\text{for}~ n \geq 2.
\end{cases}
\end{align*} 
Hence, the statement follows. 
\end{proof}

\hf \hf We encounter a similar variant of an intermediate invariant for our recursion to enumerate rational tropical curves with tangency. Let $N_d^{\mathsf{T}_0 \mathsf{T}_0}$ be the number of rational tropical degree $d$ curves intersecting $\mathcal{L}^l$ transversally at any two distinct points and passing through $3d-1$ points in general position.
Using computations similar to 
Example \ref{trop_mult_eg} and Lemma \ref{ndt0}, the number 
$N_d^{\mathsf{T}_0 \mathsf{T}_0}$ is given by
\begin{align}
\label{T0T0_num}
  N_d^{\mathsf{T}_0 \mathsf{T}_0} = l d(l d-1) ~n_d.  
\end{align}
This is because 
the two marked points get evaluated to two distinct points of $\mathcal{L}^l$ in ${l d \choose 2}~n_d$ ways with $2!$ choices to permute these points, making them distinguishable. 
We now use these numbers $N_d^{\mathsf{T}_0}(n)$ and $N_d^{\mathsf{T}_0 \mathsf{T}_0}$ for 
 proving our recursion formula for tangency.

\section{Kontsevich type formula for tangency}\label{KoTa}
This section establishes a recursion formula to enumerate rational plane tropical curves tangent to a general tropical curve of degree $l$ in the plane satisfying point constraints.  We do this by computing degree of a suitably defined tropical map at two distinct generically chosen points inside $\mathbb{R}^{2n-4}\times \M_{4}$. Recall the treatment considered in \cite{GathTropical}. The subsets $\mathcal{M}_n^{\alpha}$ of $\mathcal{M}_n$ and $\M_{\Delta, n}^{\alpha}$ of $\M_{\Delta, n}$ correspond to the marked tropical curves of type $\alpha$ and the plane tropical curve of type $\alpha$, respectively. To obtain the multiplicity of a tropical map for a given combinatorial type $\alpha$, our primary task is to find its coordinates inside $\M_{\Delta, n}^{\alpha}$ yielding an associated matrix. Then, the computation of the determinant of that associated matrix records the multiplicity of the morphism. Next, using the degree-multiplicity relation (as in \Cref{deg-mult}), we can compute the tropical degree of morphisms. Similar to \cite{GathTropical}, we need to compute degree of appropriate tropical maps.


\hf Recall the following map considered for Kontsevich's formula in \cite{GathTropical}:
\begin{defn}
    Let $n=3d$ and $d \geq 2$. Select $$\pi:= ev_1^1 \times ev_2^2 \times ev_3 \times \ldots \times ev_n \times \delta: \M_{\bigtriangleup, n} \lra  \RR^{2n-2} \times \M_{4},$$
    where $\pi$ describes the coordinates of the markings and a point in $\M_{4}$(described by the first four marked points). The map $\pi$ is a morphism of polyhedral complexes of pure dimension $2n-1$.
\end{defn}
Note that one of the crucial differences between a generic complex curve and a tropical plane curve is that it has an intrinsic multiplicity. Hence, any count of the tropical plane curve naturally corresponds to a computation of the total multiplicities. Therefore, the numbers $N_{\bigtriangleup}(\PL)$ for rational tropical curves of degree $\bigtriangleup$ through $\PL$ (collection of point constraints), 
is defined as
\begin{align}
    N_{\bigtriangleup}(\PL):= deg_{\PL}(\pi). 
\end{align}
Thus, by the Correspondence Theorem \cite[Theorem 1]{Mik-correspondence}, the authors in \cite{GathTropical} show that these numbers are the same as the corresponding classical numbers $n_d$ of rational degree $d$ stable maps to $\mathbb{CP}^2$ through $3d-1$ generic points.\\

\hf Now, we will define a suitable tropical morphism whose degree computation will produce our desired recursion formula for tangency. In what follows, by an abuse of notation, we continue calling $\overline{\mathsf{X}}$ (as opposed to $\overline{\mathsf{X}}^{trop}$) to be the tropical sub-moduli of $\M_{\bigtriangleup, n}$. Let us define the following crucial map: 

\begin{defn}
    \label{defn-tang}
    Let $n=3d$. We consider the map
\begin{align}
\label{map-giving-numbers}
    \tilde{\pi}:  ev_3 \times \ldots ev_n \times \delta :\M_{\bigtriangleup, n} \lra  \underbrace{\RR^{2}\times \ldots \times \RR^2}_{n-2} \times \M_{4},
\end{align}
where the image of 
the first two markings
lie on a common general tropical curve $\mathcal{L}^l$ as two distinct points (we will denote this situation by $\underbrace{a, b}$ on $\mathcal{L}^l$), i.e. the first two distinct coordinates lie on $\mathcal{L}^l$ (this encodes the intersection with $\overline{X}$). The image of the rest describes both the coordinates in $\RR^2$ of the other markings (since the other markings meet distinct general points). The point in $\M_{4}$ is defined by the first four marked points. Hence, the map $\tilde{\pi}$ describes all the coordinates, and is a morphism of polyhedral complexes of pure dimension $2n-3$.
\end{defn}
\hf Next, we will run the tropical WDVV algorithm developed by Gathmann-Markwig \cite[Page 550-560]{GathTropical} against the map $\tilde{\pi}$. As a geometric consequence (tropical analogue of Figure \ref{pic idea}), we obtain the numbers $N_d^{\mathsf{T}_1}(l)$ defined as the number of rational degree $d$ plane tropical curves that are tangent to $\mathcal{L}^l$ and passing through $3d-2$ generic points.
 Now, by Mikhalkin's correspondence theorem \cite{Mik-g=0-char-numbers}, these agree with the numbers 
 for rational degree $d$ complex plane curves tangent to a degree $l$ curve in $\mathbb{CP}^2$ and passing through $3d-2$ generic points in the complex plane. The agreement of the above is described in Section \ref{Verification}.\\ 
\hf For a chosen $\tilde{\pi}$-general point $\PL$, it is natural to prove the independence of $\PL$ while computing the degrees $deg_{\PL}(\tilde{\pi})$.
\begin{rem}
Throughout, our proofs will be presented using the convention that the combinatorial type of the tropical curve is rigid (see \cite{GathTropical} for dealing with non-rigid cases). 
\end{rem}
Then we have
\begin{prop}
\label{invariant-trop-degree}
    The degrees $deg_w(\tilde{\pi})$ do not depend on $w$ (whenever $w$ is in $\tilde{\pi}$-general position).
\end{prop}

\begin{proof}
    The proof follows from \cite[Proposition 4.4]{GathTropical}. However, there 
    is an extra tropical curve of degree $l$ considered, as opposed to the exact setting in \cite[Proposition 4.4]{GathTropical}. Therefore, for the sake of completeness, we will indicate the proof when the edges are bounded and the combinatorial type of the curve has exactly one $4$ valent vertex, with the remaining vertices being $3$ valent. It is sufficient to show that $deg({\tilde{\pi}})$ is locally constant around an image point under the map $\tilde{\pi}$ of a general plane tropical curve $C \in \M_{\Delta, n}^{\alpha}$ (cf. \cite[Proposition 4.4, Pg-551]{GathTropical}). \\
\hf Let the combinatorial type $\alpha$ have  four bounded edges $E_1, E_2, E_3, E_4$ around $V$ such that $E_1$ and $E_2$ both map to the degree $l$ curve considered. Let $\beta$ be a possible combinatorial type that has $\alpha$ in the boundary, as indicated in the following local picture:
\begin{figure}[hbt!]
\vspace*{0.1cm}
\begin{center}\includegraphics[scale = 0.7]{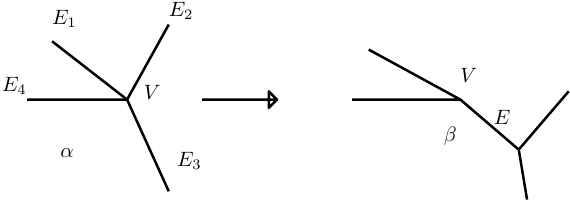}
\vspace*{-0.2cm}
\end{center}
\end{figure}
\FloatBarrier
\hf Let $l_i$ be the lengths of $E_i$ and $v_i$ be their direction vectors pointing away from $V$ for all $i = 1, \ldots, 4$, where $v_1$ is the relative direction vector measured away from $V$, assuming the distance where $x_1$ and $x_2$ maps is positive. To compute the matrix of $\tilde{\pi}$, we choose the root vertex $V$ in $\beta$, and we denote $w \in \RR^2$ to be the image of $V$. The table below encodes the relevant parts of the matrices $A$ of $\tilde{\pi}$ for the type $\beta$. The first block of the matrix consists of columns corresponding to $w$ (which is an identity matrix of size $2$) and the lengths $l_i$ for all $ i = 1, \ldots, 4$. The last column corresponds to the length of the newly added bounded edge. There are other columns corresponding to the other bounded edges except $E_i$, which are not shown here. All rows except the last row correspond to the marked points whose images are in $\RR^2$ (different types of rows reflect the location of the markings on the edges $E_i$ and how they can be reached from $V$). Note that for markings $x_i$ for $ i \geq 3 $, we use both coordinates in $\RR^2$. For $x_1$ and $x_2$, we will use only one of the coordinates of $\RR^2$, namely, we will be considering the case where $x_1$ and $x_2$ maps to the two distinct points in the $x$-axis.
The last row corresponds to the coordinate in $\M_4$, where $\star$ denotes an arbitrary value.
    \begin{table}[h!]
    \begin{tabularx}{\textwidth}{|X|X|c|}
\hline
      & $w$ \hspace*{.75cm} $l_1$ \hspace*{.75cm} $l_2$ \hspace*{.75cm} $l_3$ \hspace*{.75cm} $l_4$ & $l^{\beta}$               \\
\hline
Points behind $E_1$  & $I_2$  \hspace*{.75cm} $v_1$ \hspace*{.75cm} $0$ \hspace*{.75cm} $0$ \hspace*{.75cm} $0$ &   $0$                  \\
\hline
Points behind $E_2$  & $I_2$  \hspace*{.75cm} $0$ \hspace*{.75cm} $v_2$ \hspace*{.75cm} $0$ \hspace*{.75cm} $0$    &   $v_2 + v_3$                   \\
\hline
Points behind $E_3$  &   $I_2$  \hspace*{.75cm} $0$ \hspace*{.75cm} $0$ \hspace*{.75cm} $v_3$ \hspace*{.75cm} $0$     &   $v_2 + v_3$                   \\
\hline
Points behind $E_4$ &     $I_2$  \hspace*{.75cm} $0$ \hspace*{.75cm} $0$ \hspace*{.75cm} $0$ \hspace*{.75cm} $v_4$    &   $0$                   \\
\hline
coordinates of $\M_4$   &   $0$  \hspace*{.75cm} ~$\star$ \hspace*{.75cm} ~$\star$ \hspace*{.75cm} ~$\star$ \hspace*{.75cm} ~$\star$            &   $1$                   \\
\hline
    \end{tabularx}
\end{table}

 Next, we need to analyse the cases depending on how many edges of $C$ are contained in the subgraph $C(4)$ (a minimal connected subgraph containing $4$ unbounded edges). Consider the case when $ft_4(C)$ (for forgetful map $ft_4: \M_{\bigtriangleup, n} \lra  \M_{4}$) is the curve of type (D) in Figure \ref{type of curves in M4}, and the curve of type $\beta$ above is mapped to the type (A) curve. It follows that the matrix $A$ has all the entries of the last row as zero except the bottom right entry, which is $1$. Hence, the determinant does not depend on the entries present in the last column. Now it is straightforward to see that the matrix for all the other possible combinatorial types whose boundary contains $\alpha$ only differs at the last column (see the detailed discussion in \cite[Proposition 4.4, Pg-551]{GathTropical}). Thus, it follows that the $deg({\tilde{\pi}})$ is locally constant around $C$.  \\
 \hf Note that we assumed without loss of generality that the multiplicity of type $\beta$ is nonzero since other types are irrelevant. Thus the restriction of $\tilde{\pi}$ to $\M_{\Delta, n}^{\beta}$ is given by the invertible matrix, namely $M$. Therefore, there is at most one inverse image point in $\tilde{\pi}^{-1}(\mathcal{P})$ with coordinate $M^{-1} \mathcal{P}$. Hence, the combinatorial type $\beta$ occurs in the inverse image of the point $\mathcal{P}$ near the wall. All the other remaining cases are straightforward, and we omit the details of those.
\end{proof}

 Recall the classical WDVV procedure described in Equation \eqref{T1_WDVV}. Now if we want to obtain the equality in \eqref{T1_WDVV} using the tropical modification, all we need is to compute the degree of the map $\tilde{\pi}$ at two generically chosen points $\PL$ and $\mathcal{Q}$ with the two types of curves in $\M_4$ with a large length parameter, namely type (A) and type (B). Therefore, the classical equation \eqref{T1_WDVV} transforms into the following in the tropical setting:
\begin{align}
    \label{TWDVVT1}
 deg_{\mathcal{P}}(\tilde{\pi}) = deg_{\mathcal{Q}}(\tilde{\pi}),
\end{align}
where $\mathcal{P}, \mathcal{Q}  \in \RR^{2n-4} \times \M_4$ are two $\tilde{\pi}$ general points with $\underbrace{a, b}$ on $\mathcal{L}^l$, and having large $\M_4$ coordinates of type (A) and type (B), respectively. The above equality follows from the Proposition above.
Next, we will explain the computation of the LHS of \eqref{TWDVVT1} and only outline the computation of RHS of \eqref{TWDVVT1} due to similarity. Throughout this computation,  curves of the above types demands to acquire a contracted bounded edge, the existence of which follows from similar treatment as in \cite[Proposition 5.1]{GathTropical}.
\begin{prop}
\label{Existance-bdd-contracted-edge}
Let $d \geq 1$ and $n = 3d$. Let $\mathcal{P} \in \RR^{2n-4} \times \M_4$ be a $\tilde{\pi}$ general point with very large $\M_4$ coordinate. Then, every plane tropical curve $C$ in $\tilde{\pi}^{-1}(\mathcal{P})$ satisfying $mult_{C}(\tilde{\pi}) \neq 0$ has a bounded contracted edge.
\end{prop}
\hf Let us recall  \cite[Remark 5.2, Pg. 555]{GathTropical}. Let $C = (\Gamma, x_1, \ldots, x_n, h)$ be a degree $d$ plane tropical curve with a contracted bounded edge $E$. Then, one can split $\Gamma$ at $E$ and construct two new plane tropical curves $C_1$ and $C_2$ with mark points $x_1, \ldots, x_n$ distributed amongst them. By balancing condition, $C_i$ has degree $d_i$ with $d_1 + d_2 = d$.
Any general element in $\M_4$ with a large bounded edge corresponds to a curve of type (A), (B), or (C). Thus, geometrically, we have the following intersection theoretic situation :
\begin{lmm}
\label{multi-lemma-recursion}
    Let $\PL = (p_3, \ldots ,p_n, z) \in \RR^{2n-4} \times \M_4 $ (with $\underbrace{a, b}$ be any two distinct points on $\mathcal{L}^l$) be a point in $\tilde{\pi}$-general position and $z$ is of the type $(A)$ with a large length parameter. Then for
every tropical plane curve $C$ in $\tilde{\pi}^{-1}(\PL)$ with a nonzero multiplicity, we have the
following cases:
\begin{itemize}
    \item[(a)] The markings $x_1$ and $x_2$ are adjacent to the same vertex (with $x_1$ and $x_2$  mapping to $\underbrace{a, b}$).
    \item[(b)] Curve $C$ decomposes uniquely into two components $C_1$ and $C_2$ of degrees
$d_1$ and $d_2$ with $d_1 + d_2 = d$ such that the marked points $x_1$ and $x_2$ are on $C_1$, the points $x_3$
and $x_4$ are on $C_2$, and exactly $3d_1 - 1$ of the other points $x_5,\ldots,x_{3d}$ are on $C_1$.
\end{itemize}
\end{lmm}
\begin{proof}
    By Proposition \ref{Existance-bdd-contracted-edge}, any generic $C \in \tilde{\pi}^{-1}(\mathcal{P})$ with nonzero $\tilde{\pi}$-multiplicity have exactly one bounded contracted edge. Otherwise, we would have $mult_{C}(\tilde{\pi}) = 0$.\\
    \hf Let $E$ be the bounded contracted edge of $C$. Then, by using the tropical forgetful map, $E$ must lie in the subgraph $C(4)$ (see \cite[Lemma 5.3]{ GathTropical} for details). Consider the point $z \in \M_4$ of type (A). The markings $x_1$ and $x_2$ 
    must lie on one side of $E$ (namely $C_1$ side), and the markings $x_3$ and $x_4$ lie to the other side (namely $C_2$ side). If there are no bounded edges in $C_1$, it only consists of $E, x_1$, and $x_2$. This implies that we are in situation (a), and the evaluation conditions forces all of $C_1$ to map to $\underbrace{a, b}$.\\
    \hf Next, consider the case when there are bounded edges to both sides of $E$.  Then, the curve $C$ is reducible and we are in the situation (b). In this scenario, both $x_1$ and $x_2$ can not be on the same side of a vertex. Let $n_1$ (and respectively $n_2$) be the number of markings $x_5, \ldots, x_n$ on $C_1$ (respectively on $C_2$). We want to show $n_1 = 3d_1 -1$ and $n_2 = 3d_2 - 3$. Assume, on the contrary, that $n_1 \geq 3d_1$. Note that in this case, $2 n_1 + 2 \geq 3 d_1 + n_1 + 1 $ due to the images of the $n_1$ marked points and the first and second image points of $x_1$ and $x_2$, respectively. The  root vertex contributes $2$, giving a total of $3d_1 + n_1 + 2 - 3$ bounded edges. This leads to $mult_{C_1}(\tilde{\pi}) = 0$. Hence, we conclude that $n_1 \leq 3d_1 - 1$. The same argument as above gives us $n_2 \leq 3d_2 - 3$. Notice that the total number of points is $n_1 + n_2 = n-4$. Hence, it implies that we have the following equality $3d_1 - 1 + 3d_2 - 3 = n_1 + n_2 $. 
\end{proof}

\begin{rem}
\label{Gluing_remark}
Note the converse of the above lemma is also true. For any point $\PL= (p_3, \ldots ,p_n, z)$ in $\tilde{\pi}$-general position with  $\underbrace{a, b}$ on $\mathcal{L}^l$, and having $\M_4$ coordinate of type (A) with a large length parameter, consider two plane tropical curves $C_1$ and $C_2$ of degree $d_1$ and $d_2$ respectively with $d_1+d_2 = d$ such that the image of $C_1$ passes through $3d_1 - 1$ general points from $p_5, \ldots ,p_n$, and intersects $\mathcal{L}^l$ at $\underbrace{a, b}$. 
Further, the image of $C_2$ passes through $p_3$ and $p_4$ along with $3d_2 - 3$ other general points. Then, for each choice of points $P\in C_1$ and $Q \in C_2$ that map to the same image point in $\RR^2$, and for each choice of markings $x_1, \ldots, x_n$ on $C_1$ and $C_2$ that map to  $\underbrace{a, b}$ on $\mathcal{L}^l$ and to the general points $p_3, \ldots, p_n$, respectively, one can construct a single reducible $n$ marked curve $C \in \tilde{\pi}^{-1}(\PL)$ by gluing $C_1$ and $C_2$. The curve is $3$ valent if the components are so. 
In other words, choosing  $C_1$, $C_2$ and respective points $x_1, \ldots x_n, P, Q$, there is a unique curve in $\tilde{\pi}^{-1}(\PL)$ obtained from this data. Thus, while computing degree of $\tilde{\pi}$ for the curves occurring in situation $(b)$ of \Cref{multi-lemma-recursion}, we can just sum over all choices $C_1, C_2, x_1, \ldots, x_n, P, Q$ as discussed above.
\end{rem}

\hf Now we are left with computing the degree of $\tilde{\pi}$. Before that we must compute the multiplicities of $\tilde{\pi}$ for all possible curves in $\tilde{\pi}^{-1}(\PL)$.

\begin{prop}
\label{multi-of-the-map-pi}
    With the notations used in Lemma $\ref{multi-lemma-recursion}$ and Proposition $\ref{invariant-trop-degree}$, let $C$ be a point in $\tilde{\pi}^{-1}(\PL)$. Then
    \begin{itemize}
        \item[(i)] If $C$ is considered as in situation $(a)$ of Lemma \ref{multi-lemma-recursion}, its $\tilde{\pi}$-multiplicity is $mult_{C^{\prime}}(ev)$, where $C^{\prime}$ is
obtained from $C$ by forgetting the first two markings, and $ev$ is the product of evaluations at the $3d-2$ points $x_3,\ldots,x_{3d}$. Thus the multiplicity counts the number of degree $d$ tropical rational curves tangent to $\mathcal{L}^l$ (two distinct points in $\mathcal{L}^l$  come together).
\item[(ii)] If $C$ is considered as in situation $(b)$ of Lemma \ref{multi-lemma-recursion}, its $\tilde{\pi}$-multiplicity is 
\begin{align}
\label{Evalution-mult}
mult_{C}(\tilde{\pi}) = mult_{C_1}(ev) ~mult_{C_2}(ev)~  \left(C_1 \cdot C_2 \right)_{P=Q}, 
\end{align}
 \end{itemize}
 where $\left(C_1 \cdot C_2 \right)_P$ denotes the intersection multiplicity of the tropical curves $C_1$ and $C_2$ at $P \in C_1 \cap C_2$. 
\end{prop}

\begin{proof}
  We need to compute the determinant of the associated matrix for the map $\tilde{\pi}$. Note that the length of the bounded contracted edge is irrelevant for the evaluations in both the cases since it contributes $1$ to the $\M_4$ coordinate of $\tilde{\pi}$. Thus it is sufficient to compute the determinant by dropping the row coming from the $\M_4$ coordinate and the column corresponding to the bounded contracted edge. Also, the multiplicity calculation does not depend on the choice of the root vertex.\\
    \hf First case (i) follows by treating the images for $x_1$ and $x_2$ as the image of a single
    marked point as opposed to computing the first evaluation at $x_1$ and then the second evaluation at $x_2$. This is because the tripod consists of a contracted edge and the adjacent markings $x_1$ and $x_2$ getting contracted to a point along $\mathcal{L}^l$. Thus, replacing the combinatorial type with a contracting end along $\mathcal{L}^l$  gives the computation for rational degree $d$ tropical curves tangent to $\mathcal{L}^l$ passing through $3d-2$ general points.\\  
 \hf Finally, for the case (ii), consider the bounded contracted edge $E$, at which by splitting the curve $C$, we have two components $C_1$ and $C_2$ with the root vertex in one of them, namely on the $C_1$ side. 
 Denote the length coordinates by $l_i = l(E_i)$. Next, the image of the markings $x_1$ and $x_2$ lie on $\mathcal{L}^l$.
 Let $E_1$ and $E_2$ be two bounded adjacent edges with common direction vector $v = (v_1, v_2)$
    as depicted below.
\begin{figure}[hbt!]
\vspace*{0.1cm}
\begin{center}\includegraphics[scale = 0.7]{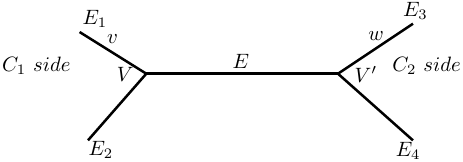}
\vspace*{-0.2cm}
\end{center}
\end{figure}
\FloatBarrier    
  Assuming the root vertex is on the side of $E_1$, the entries of the matrix $\tilde{\pi}$ corresponding to the above local picture is given by\\
     
      \begin{table}[h!]
    \begin{tabularx}{\textwidth}{|X|c|X||c||X|}
\hline
     &$h(V)$ & Lengths in $C_1$\hspace*{.45cm} $l_1$\hspace*{.25cm} $l_2$ & $l$ & Lengths in $C_2$\hspace*{.45cm}$l_3$\hspace*{.25cm} $l_4$                 \\
\hline
Points behind $E_1$ & $I_2$  & \hspace*{0.3cm} $\star$ \hspace*{2.0cm} $v$ \hspace*{.25cm} $0$   & $0$  &  \hspace*{0.3cm} $0$ \hspace*{2.0cm} $0$ \hspace*{.25cm} $0$   \\
\hline
Points behind $E_2$  & $I_2$   & \hspace*{0.3cm} $\star$ \hspace*{2.0cm} $0$ \hspace*{.25cm} $-v$   & $0$  &  \hspace*{0.3cm} $0$ \hspace*{2.0cm} $0$ \hspace*{.25cm} $0$ \\
\hline
Points behind $E_3$  &   $I_2$   &  \hspace*{0.3cm} $0$ \hspace*{2.0cm} $0$ \hspace*{.25cm} $0$   & $0$  &  \hspace*{0.3cm} $\star$ \hspace*{2.0cm} $w$ \hspace*{.25cm} $0$   \\
\hline
Points behind $E_4$ &     $I_2$   &   \hspace*{0.3cm} $0$ \hspace*{2.0cm} $0$ \hspace*{.25cm} $0$   & $0$  &  \hspace*{0.3cm} $\star$ \hspace*{2.0cm} $0$ \hspace*{.25cm} $-w$      \\
\hline
coordinates of $\M_4$   &   $0$            &  \hspace*{0.3cm} $0$ \hspace*{2.0cm} $0$ \hspace*{.25cm} $0$   & $1$  &  \hspace*{0.3cm} $0$ \hspace*{2.0cm} $0$ \hspace*{.25cm} $0$    \\
\hline
    \end{tabularx}
\label{table_exp}
\end{table}
\FloatBarrier  
(see \cite{Renzo_notes} for more details). Note that we have displayed only a part of the $6d -3 \times 6d -3$ matrix. Following the same techniques adopted in \cite{Renzo_notes}, we get the following similar block matrices $A_1$ and $A_2$, where $n_1, n_2$ are as in \Cref{multi-lemma-recursion} and $\star$ denotes arbitrary entry. Then using column operation, we get the following block matrices
\begin{enumerate}
\item[(I)] The matrix $A_1$ is obtained from the coordinates coming from the moduli point associated with $C_1 = (\Gamma_1, x_1, \ldots, x_{n_1}, h_1)$ with the markings $x_1$ and $x_2$ mapping to $\underbrace{a, b}$ in $\mathcal{L}^l$ through $h_1$. Consider the block consisting of the first four rows. It can be shown that this is similar to the block matrix $A_1$ of size $6d_1-1 \times 6d_1-1$, which computes the multiplicity of $\tilde{\pi}$ for $C_1$:\\
\begin{center}
$A_1 = $
    \begin{tabular}{c||l}
\hline
    & $V$ \hspace*{.65cm} Lengths in $C_1$\hspace*{.65cm} $l_1$ \hspace*{.65cm} $\tilde{l}_2$ \hspace*{.65cm} $\tilde{l}_3$ \hspace*{.65cm} $\tilde{l}_4$  \\
\hline
   Points behind $E_1$ & $I_2$  \hspace*{1.15cm} $\star$ \hspace*{2cm} $v$ \hspace*{.71cm} $0$ \hspace*{.71cm} $0$ \hspace*{.63cm} $0$                  \\
\hline
 Points behind $E_2$  & $I_2$  \hspace*{1.15cm} $\star$ \hspace*{2cm} $0$ \hspace*{.75cm} $0$ \hspace*{.65cm} $0$ \hspace*{.65cm} $0$     \\
\hline
    \end{tabular}
    \end{center}
  \vspace*{.4cm} 
where $\tilde{l}_2,\tilde{l}_3$ and $\tilde{l}_4$ are obtained after suitable column operation. This implies that the determinant of the matrix gives us the following multiplicity
\begin{align}
\label{first_block_mult}
det(A_1) = mult_{C_1}(\tilde{\pi}) 
\end{align}
\item[(II)] Similar to the above, consider the lower block of the matrix. It is similar to the matrix $A_2$ of size $3d_2-2 \times 3d_2-2$, which compute the multiplicity of the $\tilde{\pi}$ for $C_2 = (\Gamma_2, x_{n_1+1}, \ldots, x_{n}, h_2)$:  the matrix consists of $2 \times 2$ block of the vectors $v, w$ and a lower identity block. Then\\
\begin{center}
$A_2 = $
    \begin{tabular}{c||l|c}
\hline
    & $V^{\prime}$ \hspace*{.75cm} $\tilde{l}_1$ \hspace*{.75cm} $\tilde{l}_2$ \hspace*{.75cm} $l_3$ \hspace*{.75cm} $\tilde{l}_4$ & Lengths in $C_2$                \\
\hline
   Points behind $E_3$ & $I_2$  \hspace*{.75cm} $0$ \hspace*{.75cm} $0$ \hspace*{.75cm} $v$ \hspace*{1.1cm} $0$ &   $\star$                  \\
\hline
 Points behind $E_4$  & $I_2$  \hspace*{.75cm} $0$ \hspace*{.75cm} $0$ \hspace*{.75cm} $w$ \hspace*{.75cm} $-w$    &   $\star$                   \\
\hline
    \end{tabular}
    \end{center}
 \vspace*{.5cm}  
where $\tilde{l}_1,\tilde{l}_2$ and $\tilde{l}_4$ are obtained after suitable column operation. Using the idea applied in the proof of \cite[Proposition 3.8]{GathTropical}, we get\[    |det(A_2) | = \left( C_1 \cdot C_2 \right)_{P = Q}   \cdot mult_{C_2}(ev).   \]
\end{enumerate}   

Gathering the above two multiplicities, we completed the formula for multiplicity \eqref{Evalution-mult}.  
    
\end{proof}

\hf Next, we require analogous statements to \Cref{multi-lemma-recursion}, \Cref{Gluing_remark}, and \Cref{multi-of-the-map-pi} for the (B) type curves in $\M_4$  to compute the RHS of the Equation \eqref{TWDVVT1}. Note that for the analogous lemma to \Cref{multi-lemma-recursion}, the situation (a) in \Cref{multi-lemma-recursion} does not contribute to the recursion since the two distinct markings  map to the intersection of a generic line and a point. Note that this intersection is empty. Thus only situation (b) in \Cref{multi-lemma-recursion} is relevant. Consider another $\tilde{\pi}$-general point $\mathcal{Q} = (q_3, \ldots, q_n, z) \in \RR^{2n-4} \times \M_4$ with $\underbrace{a^{'}, b^{'}}$ on $\mathcal{L}^l$. However, in this case $a^{'}$ lies in the intersection of the image of $C_1$ and $\mathcal{L}^l$ and $b^{'}$ lies in the intersection of the image of $C_2$ and $\mathcal{L}^l$. Hence, in this case, the analogous proposition to \Cref{multi-of-the-map-pi} is as follows:
\begin{prop}
\label{RHSmulti-of-the-map-pi}
    With the above notations,
    let $C$ be a point in $\tilde{\pi}^{-1}(\mathcal{Q})$. 
If $C$ is considered as in situation $(b)$ of the lemma analogous to \Cref{multi-lemma-recursion}, and for the type (B) curve in $\M_4$ similar to Lemma \ref{multi-lemma-recursion}, then its $\tilde{\pi}$-multiplicity is 
\begin{align}
\label{RHSEvalution-mult}
mult_{C}(\tilde{\pi}) = mult_{C_1}(ev) ~mult_{C_2}(ev)~  \left(C_1 \cdot C_2 \right)_{P=Q} \cdot \left( C_1 \cdot \mathcal{L}^l\right)_{x_1} \cdot \left( C_2 \cdot \mathcal{L}^l\right)_{x_2},
\end{align}
 where $\left(C_1 \cdot C_2 \right)_P$ denotes the intersection multiplicity of the tropical curves $C_1$ and $C_2$ at $P \in C_1 \cap C_2$. The first and second coordinates of the direction vectors of $C_1$ and $C_2$ are denoted by $\left(C_i \cdot \mathcal{L}^l \right)_{x_i}$ for all $i = 1, 2$.
\end{prop}
\begin{proof}

The proof is similar to the proof of \Cref{multi-of-the-map-pi}.

\end{proof}

\hf We now prove our main theorem giving a recursion for the number $N_d^{\mathsf{T}_1}(l)$ of rational degree $d$ tropical plane curves that are tangent to $\mathcal{L}^l$ and passing through $3d-2$ generic points.
\begin{thm}
\label{final theo}
    Let $l, d \geq 1$. The numbers $N_{d}^{\mathsf{T}_1}(l)$ satisfy the following
recursion formula
\begin{align}
\label{final recursion for T1}
     N_{d}^{\mathsf{T}_1}(l)  =  \sum_{\substack{  d_1 + d_2 = d\\ d_1, d_2 > 0}} \binom{3d-4}{3d_1-2} l d_1~ l d_2~  d_1d_2 ~N_{d_1}^{\big(\mathsf{T}_0\big)_{\textnormal{pt}}} ~N_{d_2}^{\big(\mathsf{T}_0\big)_{\textnormal{pt}}} - \binom{3d-4}{3d_2-3}  d_1 d_2 ~N_{d_1}^{\mathsf{T}_0 \mathsf{T}_0} ~ n_{d_2}.
\end{align}

\end{thm}

Note that the numbers $n_d$ for plane tropical rational degree $d$ curves play the role of the base case for the above recursion. The number $ N_d^{\mathsf{T}_{1}}(1)$ for $d=1$ is $0$. This is because of a simple geometric reason: an irreducible line can not be tangent to a fixed tropical line in the plane.

\begin{proof}
    We prove the theorem by computing the tropical degree of $\tilde{\pi}$ defined in \eqref{map-giving-numbers} at two distinct general points. Let $\PL = (p_3,...,p_n, z) \in \RR^{2n-4} \times \M_4$ with $\underbrace{a, b}$ on $\mathcal{L}^l$ as two distinct points. The curve in $\M_4$ considered is of type $(A)$ (see Figure \ref{type of curves in M4}).
    Next, we count the points of $\tilde{\pi}^{-1}(\PL)$ with its $\tilde{\pi}$-multiplicity. Note that for curves considered in situation $(a)$ of Lemma \ref{multi-lemma-recursion}, we can enumerate curves of degree $d$ through $3d -2$ points with their
$\tilde{\pi}$-multiplicity using Proposition \ref{multi-of-the-map-pi}. This gives us the number $ N_d^{\mathsf{T}_{1}}(l)$.

\hf Next, for curves considered in situation $(b)$ of Lemma \ref{multi-lemma-recursion}, we need to count all possible $(C_1, C_2, x_1, \ldots x_n, P, Q)$, where
\begin{enumerate}
\item[$\bullet$] $C_1, C_2$ have degrees $d_1, d_2$, respectively such that $d_1 + d_2 = d$.
\item[$\bullet$] $x_1, x_2$ markings are on $C_1$ and maps to $\mathcal{L}^l$ at two distinct points, respectively.
\item[$\bullet$] $x_3, x_4$ markings are on $C_2$ that maps to $p_3, p_4$, respectively.
\item[$\bullet$] The markings $x_5, \ldots, x_n$ are mapped to $p_5, \ldots, p_n$, out of which $3d_1 -1$ of them lie on $C_1$ and $3d_2 -3$ of them lie on $C_2$.
\item[$\bullet$] The points $P\in C_1$ and $Q\in C_2$ (as in \Cref{Gluing_remark}) are mapped to the nodal point in $\RR^2$.
\end{enumerate}
 Now we have ${3d-4\choose 3d_2-3}$ choices to split up the points $x_5, \ldots, x_n$. Let us consider the two following sets:\\
\begin{itemize}
\item[$\bullet$] $ X_1:= \{(\Gamma_1, h_1) \in \M_{ d_1, 3d_1+1} ~|~  h_1(x_1),~ h_1(x_2) \in \mathcal{L}^l, ~h_1(x_i) = p_i, ~\forall i = 5,\ldots, 3d_1+1 \} $,\\
\item[$\bullet$] $ X_2:= \{(\Gamma_2, h_2) \in \M_{ d_2, 3d_2-1} ~|~ h_1(x_3)=p_3,~ h_1(x_4)=p_4,~~  ~h_2(x_i) = p_{i}, ~\forall i = 5,\ldots, 3d_2-1 \}. $
\end{itemize}
From earlier discussion and using \eqref{T0T0_num} and \eqref{TKon_num}, we compute the weighted sums
\begin{align}
\sum_{X_1} mult_{C_1}(\tilde{\pi}) = N_{d_1}^{\mathsf{T}_{0}  \mathsf{T}_{0}},
\end{align}
and 
\begin{align}
\sum_{X_2} mult_{C_2}(\tilde{\pi}) = n_{d_2}.
\end{align}
If we want to count $\tilde{\pi}$-multiplicity  for $C_1$ and $C_2$ as in situation $(b)$ of Lemma \ref{multi-lemma-recursion}, then fixing the curve components of degree $d_1$ and $d_2$, we get $N_{d_1}^{\mathsf{T}_{0}  \mathsf{T}_{0}} ~ n_{d_2}$ 
choices (see Proposition \ref{multi-of-the-map-pi}). By tropical B\'{e}zout’s theorem (see \cite[Theorem 1.3.2]{Strum-book}), the component $C_1$ of degree $d_1$ intersects the components $C_2$ of degree $d_2$ at the nodal point i.e., when $P = Q$ with the local multiplicity $\left(C_1 \cdot C_2 \right)_{P=Q}$ and hence there are $d_1 d_2$ choices for that.
Combining all the above, we have the degree of $\tilde{\pi}$ at $\PL$ as
$$ deg_{\PL}(\tilde{\pi}) = N_{d}^{\mathsf{T}_1} + \sum_{\substack{  d_1 + d_2 = d\\ d_1, d_2 > 0}} \binom{3d-4}{3d_2-3}~  d_1 d_2 ~N_{d_1}^{\mathsf{T}_{0}  \mathsf{T}_{0}} ~ n_{d_2}.$$ 
Next, we repeat the same process for another point $\mathcal{Q}$ with $\mathcal{M}_4$ coordinate of type $(B)$. Then, in a similar manner, using \Cref{RHSmulti-of-the-map-pi} and \Cref{T0_num}, we get the degree of $\tilde{\pi}$ at $\mathcal{Q}$ as
$$ deg_{\mathcal{Q}}(\tilde{\pi}) = \sum_{\substack{  d_1 + d_2 = d\\ d_1, d_2 > 0}} \binom{3d-4}{3d_1-2}~ l d_1~ l d_2~ d_1 d_2~ N_{d_1}^{\left( \mathsf{T}_0 \right)_{\textnormal{pt}}}~ N_{d_2}^{\left( \mathsf{T}_0 \right)_{\textnormal{pt}}}.$$
Next, by equating the two expressions above, 
we get our desired recursion formula \eqref{final recursion for T1}. This completes the proof.

\end{proof}

\section{Verification with existing results}
\label{Verification}
In this section, we verify the numbers obtained from our recursion formula \eqref{final recursion for T1} with various existing results. 
In the classical sense i.e., in $\mathbb{CP}^2$, using complex and symplectic geometric techniques, counting rational curves with tangency have already been studied. Amongst others, this includes \cite{Ernstroem1998a}, \cite{rv1}, \cite[Section 4]{Pan@Qdiv}, where they considered the case where the complex plane rational curves are tangent to a line at some point (i.e., the point is free to move in the line). In this case, the classical characteristic numbers are 
\begin{center}
\vspace{.2cm}
\renewcommand{\arraystretch}{1.5}
\begin{tabular}{|c|c|c|c|c|c|c|c|c|} 
\hline
$d$ &$1$ &$2$ & $3$ & $4$ & $5$ & $6$ & $7$ & $8$ \\
\hline 
$N_d^{\mathsf{T}_{1} }(1)$ 
& $0$ & $2$ & $36$ & $2184$ & $335792$ & $106976160$ & $61739450304$ & $ 58749399019136 $ \\ 
\hline
\end{tabular}
\vspace{.2cm}
\end{center}
These are in agreement with the numbers obtained by our recursion \eqref{final recursion for T1} when $l=1$. We have implemented our recursion formula \eqref{final recursion for T1} in Mathematica software, available in
\begin{center}
\url{https://github.com/AnantaSpace/Tropical-tangency-via-tropical-WDVV.git}\\
\end{center}
\hf Nevertheless, from our recursion \eqref{final recursion for T1}, one can very easily extract the numbers $N_d^{\mathsf{T}_{1} }(l)$ for all $l\geq 1$ (using the above code). Classically, Gathmann has enumerated the number of complex degree $d$ rational curves that are tangent to a degree $l$ curve in $\mathbb{CP}^2$ satisfying point conditions. We will tabulate various cases for different values of $l$, which are extensively verified by the computer program ``GROWI'' developed by Gathmann (see \cite{GROWI}). 
\begin{center}
\renewcommand{\arraystretch}{1.5}
\vspace{.2cm}
\begin{tabular}{|c|c|c|c|c|c|c|c|c|} 
\hline
$d$ &$1$ &$2$ & $3$ & $4$ & $5$ & $6$ & $7$ & $8$ \\
\hline 
$N_d^{\mathsf{T}_{1} }(2)$ 
& $0$ & $6$ & $96$ & $5608$ & $846192$ & $266578272$ & $152712516992$ & $ 144550300093056 $ \\ 
\hline
$N_d^{\mathsf{T}_{1} }(3)$ 
& $0$ & $12$ & $180$ & $10272$ & $1531200$ & $478806336$ & $272919200064$ & $ 257402703221760 $ \\ 
\hline 
$N_d^{\mathsf{T}_{1} }(4)$ 
& $0$ & $20$ & $288$ & $16176$ & $2390816$ & $743660352$ & $422359499520$ & $ 397306608405248$ \\ 
\hline 
$N_d^{\mathsf{T}_{1} }(5)$ 
& $0$ & $30$ & $420$ & $23320$ & $3425040$ & $1061140320$ & $601033415360$ & $564262015643520 $ \\ 
\hline
\end{tabular}
\vspace{.2cm}
\end{center}

\section*{Acknowledgement}
We thank Bivas Khan for several useful discussions. 

\bibliographystyle{siam}

\end{document}